\documentclass[a4paper,11pt]{article}

\usepackage{amsmath,amssymb,latexsym}
\usepackage{graphicx}
\usepackage{epsfig}
\usepackage{longtable}
\usepackage{amsthm}
\usepackage{color}

\newtheorem{theorem}{Theorem}

\newtheorem{proposition}[theorem]{Proposition}

\newtheorem{remark}{Remark}

\begin{document}

\begin{center}
{\Large\bf Best-possible bounds on the set of copulas with a given value of Gini's gamma}
\end{center}

\begin{center}
{\large\bf Manuel \'Ubeda-Flores}
\end{center}

\begin{center}
{\it Department of Mathematics, University of Almer\'{\i}a, 04120 Almer\'{\i}a, Spain.\\
{\rm mubeda@ual.es}}
\end{center}
\smallskip

\begin{center}{\bf Abstract}\end{center}

\noindent In this note, pointwise best-possible (lower and upper) bounds on the set of copulas with a given value of the Gini's gamma coefficient are established. It is shown that, unlike the best-possible bounds on the set of copulas with a given value of other known measures such as Kendall's tau, Spearman's rho or Blomqvist's beta, the bounds found are not necessarily copulas, but proper quasi-copulas.
\bigskip

\medskip


\section{Introduction}

Quasi-copulas were introduced in the field of probability in order to characterize operations on distribution functions that can or cannot be derived from operations on random variables defined on the same probability space (see \cite{Alsina93}), and were characterized in \cite{Genest99}. In the last years  these functions have attracted an increasing interest by researchers in some topics of fuzzy sets theory, such as preference modeling, similarities and fuzzy logics. For a survey on quasi-copulas, see \cite{Sempi2017}.

Copulas, probability distribution functions with uniform margins on $[0,1]$, are a subclass of quasi-copulas. The importance of copulas in probability and statistics comes from {\it Sklar's theorem} \cite{Sklar59}, which states that the joint distribution $H$ of a pair of random variables $(X,Y)$---defined on the same probability space $(\Omega,{\cal{F}},\mathbb{P})$---and the corresponding (univariate) marginal distributions $F$ and $G$ are linked by a copula $C$ in the following manner:
$$H(x,y)=C\left(F(x),G(y)\right)\,\, {\rm for}\,\,{\rm all}\,\, \left(x,y\right)\in[-\infty,\infty]^2.$$
If $F$ and $G$ are continuous, then the copula is unique. For a review on copulas, we refer to the monographs \cite{Durante2016book,Nelsen2006}.

The fundamental best-possible bounds inequality for the set of (quasi-)copulas is given by the Fr\'echet-Hoeffding bounds, i.e., for any quasi-copula $Q$ we have
\begin{equation}\label{Frechet}
W(u,v):=(0\vee (u+v-1))\le Q(u,v)\le (u\wedge v)=:M(u,v)
\end{equation}
for all $(u,v)\in[0,1]^2$, where $c\vee d:=\max(c,d)$ and $c\wedge d:=\min(c,d)$ for any two real numbers $c$ and $d$. Furthermore, the bounds $W$ and $M$ are themselves copulas.

A procedure for finding pointwise best-possible bounds on sets of copulas and a given value of the population version of a measure of association, such as Kendall's tau, Spearman's rho and the population version of the medial correlation coefficient (or Blomqvist's beta) is illustrated in \cite{NelUb04,Nel01}. The bounds attained are evaluated, with the result that all the bounds are copulas. In the case of the Spearman's footrule coefficient, the lower bound is a copula, but the upper bound can be a copula or a proper quasi-copula (see \cite{Beliakov2022,Kokol2021}).

In this note we focus on the Gini's gamma coefficient and establish the best-possible (lower and upper) bounds on the set of copulas with a given value of Gini's gamma. The resulting bounds are not always copulas, but proper quasi-copulas. The construction of the bounds follows the same structure as for the Spearman's footrule bounds in \cite{Beliakov2022} (for a different idea, where the bounds are expressed in an equivalent form, see \cite{Kokol2021}).

After some preliminaries concerning (quasi-)copulas (Section \ref{prelim}), we present the main results in Section \ref{mainsec}, where we find the best-possible bounds on the set of copulas with a given value of the Gini's gamma coefficient and provide some salient properties. Section \ref{conclus} is devoted to conclusions.

\section{Preliminaries}\label{prelim}

A (bivariate) {\it copula} is a function $C\colon [0,1]^{2}\longrightarrow [0,1]$ which satisfies
\begin{itemize}
\item[(C1)] the boundary conditions $C(t,0)=C(0,t)=0$ and $C(t,1)=C(1,t)=t$ for all $t$ in $[0,1]$, and
\item[(C2)] the {\it 2-increasing property}, i.e., $V_{C}([u_{1},u_{2}]\times[v_{1},v_{2}])=C(u_{2},v_{2})-C(u_{2},v_{1})-C(u_{1},v_{2})+C(u_{1},v_{1})\ge 0$ for all $u_{1},u_{2},v_{1},v_{2}$ in $[0,1]$ such that $u_{1}\le u_{2}$ and $v_{1}\le v_{2}$.
\end{itemize}

$V_C(R)$ is usually called as the $C$-{\it volume} of the rectangle $R$; and in the sequel we also consider the $C$-volume of a rectangle for real-valued functions on $[0,1]^n$ which may not be copulas.

Let $\mathcal{C}$ denote the set of all copulas.

A (bivariate) {\it quasi-copula} is a function $Q\colon[0,1]^{2}\longrightarrow [0,1]$ which satisfies condition (C1) of copulas, but in place of (C2), the weaker conditions
\begin{itemize}
\item[(Q1)] $Q$ is non-decreasing in each variable, and
\item[(Q2)] the {\it Lipschitz} condition $|Q(u_{1},v_{1})-Q(u_{2},v_{2})|\le |u_{1}-u_{2}|+|v_{1}-v_{2}|$ for all $u_{1},v_{1},u_{2},v_{2}$ in $[0,1]^{2}$.
\end{itemize}

While every copula is a quasi-copula, there exist {\it proper} quasi-copulas, i.e., quasi-copulas which are not copulas.

One of the most important occurrences of quasi-copulas in statistics is due to the following observation (\cite{NelUb05,Nelsen04}): Every set $\cal{S}$ of (quasi-)copulas has the smallest upper bound and the greatest lower bound in the set of quasi-copulas (in the sense of pointwisely ordered functions). These bounds do not necessarily belong to the set $\cal{S}$, nor they are necessarily copulas if the set consists of copulas only.

In \cite{Tankov} best-possible bounds on the set of quasi-copulas that coincide on a given compact subset $S$ of $[0,1]^2$ are established and the author investigates sufficient conditions on $S$ such that these bounds are also the best-possible bounds on the set of copulas that coincide on $S$.

The bounds in \eqref{Frechet} can often be narrowed when we possess additional information about the copula $C$. In \cite{Nelsen2006} the best-possible bounds on a set of copulas when a value at a single point is known are provided (we recall this result for our purposes, and in which $x^+:=0\vee x$ for any real number $x$).

\begin{proposition} \label{prop1}
Let $C$ be a copula, and suppose $C(a,b)=\theta$, where $(a,b)\in [0,1]^2$ and $W(a,b)\le \theta\le M(a,b)$. Then $\underline{C}_{(a,b),\theta}(u,v)\le C(u,v)\le \overline{C}_{(a,b),\theta}(u,v)$ for all $(u,v)\in[0,1]^2$, where
\begin{equation*}\label{Lower}
\underline{C}_{(a,b),\theta}(u,v)=\max\left(0,u+v-1,\theta-(a-u)^+-(b-v)^+\right)
\end{equation*}
and
\begin{equation*}\label{Upper}
\overline{C}_{(a,b),\theta}(u,v)=\min\left(u,v,\theta+(u-a)^++(v-b)^+\right).
\end{equation*}
Since $\underline{C}_{(a,b),\theta}(a,b)=\overline{C}_{(a,b),\theta}(a,b)=\theta$, the bounds, which are copulas, are best-possible.
\end{proposition}

\section{Best-possible bounds when a given value of Gini's gamma is known}\label{mainsec}

Let $(R_1, S_1),\ldots,(R_n, S_n)$ be ranks associated with a random sample $(X_1, Y_1),\ldots,(X_n, Y_n)$ from some
continuous bivariate distribution $H(x, y) = {\mathbb{P}}(X \le x, Y \le y)$. The Italian statistician Corrado Gini \cite{Gini} introduced the {\it indice di cograduazione semplice}---also known as the {\it Gini's rank association coefficient}---in the following manner:
$$\gamma_n=\frac{1}{\lfloor n^2/2\rfloor}\sum_{i=1}^{n}\left\{|n+1-R_i-S_i|-|R_i-S_i|\right\},$$
where $\lfloor x\rfloor$ denotes the integer part of $x>0$. In terms of the copula $C$ associated with the continuous random vector $(X,Y)$, the {\it Gini's gamma}---as it is also known---, denoted by $\gamma(C)$, can be expressed as
$$\gamma(C)=4\int_{0}^{1}\left[C(u,u)+C(u,1-u)\right]\,du -2$$
(see \cite{NelsenGini,Nelsen2006}). We note that, for any copula $C$, we have $-1=\gamma(W)\le\gamma(C)\le \gamma(M)=1$. For an overview of this coefficient, we refer to \cite{GenestGini}.

For any $t\in[-1,1]$, let ${\bf G}_t$ denote the set of copulas with a common value $t$ of Gini's gamma, i.e.,
$${\bf G}_t = \left\{C\in{\cal{C}}\colon \gamma(C)=t\right\}.$$
Let $\underline{G}_t$ and $\overline{G}_t$ denote, respectively, the pointwise infimum and supremum of ${\bf G}_t$, i.e., for each $(u,v)\in[0,1]^2$,
\begin{equation}\label{InfsupG}
\underline{G}_t(u,v) =\inf\left\{C(u,v)\colon C\in{\bf G}_t\right\}\,\,{\rm and}\,\, \overline{G}_t(u,v) = \sup\left\{C(u,v)\colon C\in{\bf G}_t\right\}.
\end{equation}

A relationship between $\underline{G}_t$ and $\overline{G}_t$ is provided in the next result.

\begin{proposition}\label{relation}
Let $\underline{G}_t$ and $\overline{G}_t$ be the pointwise infimum and supremum \eqref{InfsupG} of ${\bf G}_t$. Then we have
$$\underline{G}_t(u,v)=v-\overline{G}_{-t}(1-u,v),$$
or equivalently,
$$\underline{G}_t(u,v)=u-\overline{G}_{-t}(u,1-v)$$
for all $(u,v)\in[0,1]^2$.
\end{proposition}

\begin{proof}We prove $\underline{G}_t(u,v)=v-\overline{G}_{-t}(1-u,v)$ for all $(u,v)\in[0,1]^2$ (the proof for the other equality is similar and we omit it). Let $C$ be a copula, and let $\widetilde{C}$ be the copula given by $\widetilde{C}(u,v)=v-C(1-u,v)$ for every $(u,v)$ in $[0,1]^2$, i.e., if $(X,Y)$ is the random pair whose associated copula is $C$ then $(1-X,Y)$ is the pair associated with the copula $\widetilde{C}$. Then
\begin{eqnarray*}
\gamma\left(\widetilde{C}\right)&=&4\int_{0}^{1}\left[\widetilde{C}(u,u)+\widetilde{C}(u,1-u)\right]du-2\\
&=&4\int_{0}^{1}\left[u-C(1-u,u)+1-u-C(1-u,1-u)\right]du-2\\
&=&2-4\int_{0}^{1}[C(1-u,1-u)+C(1-u,u)] du\\
&=&2-4\int_{0}^{1}[C(z,z)+C(z,1-z)] dz\\
&=&-\gamma(C).
\end{eqnarray*}
Now, using this equality, note that we have
\begin{eqnarray*}
\underline{G}_t(u,v)&=&\inf\left\{C(u,v)\,:\, C\in\mathcal{C},\, \gamma(C)=t\right\}\\
&=&\inf\left\{C(u,v)\,:\,C\in\mathcal{C},\, \gamma\left(\widetilde{C}\right)=-t\right\}\\
&=&\inf\left\{\widetilde{C}(u,v)\,:\, C\in\mathcal{C},\, \gamma\left(C\right)=-t\right\}\\
&=&\inf\left\{v-C(1-u,v)\,:\, C\in\mathcal{C},\, \gamma\left(\widetilde{C}\right)=-t\right\}\\
&=&v-\sup\left\{C(1-u,v)\,:\, C\in\mathcal{C},\, \gamma\left(\widetilde{C}\right)=-t\right\}\\
&=&v-\overline{G}_{-t}(1-u,v),
\end{eqnarray*}
and this completes the proof.
\end{proof}

The next result provides explicit expressions for the bounds $\underline{G}_t$ and $\overline{G}_t$.

\begin{theorem}\label{boundsgini}Let $\underline{G}_t$ and $\overline{G}_t$ be the pointwise infimum and supremum \eqref{InfsupG} of ${\bf G}_t$ for $t\in [-1,1]$. Then, for any $(u,v)\in [0,1]^2$, we have

\begin{equation}\label{Gsup}
\overline{G}_t(u,v)=\min\left(u,v,\max\left(\chi_{\theta_1}\theta_1,\chi_{\theta_2}\theta_2,\chi_{\theta_3}\theta_3,\chi_{\theta_4}\theta_4,\chi_{\theta_5}\theta_5\right)\right)
\end{equation}
and
\begin{equation}\label{Ginf}
\underline{G}_t(u,v)=v-\overline{G}_{-t}(1-u,v),
\end{equation}
where
\begin{eqnarray}\label{theta1}
\theta_1\!\!\!&=&\!\!\!\frac{u+v-1+\sqrt{(u+v-1)^2+t+1}}{2},\\\label{theta2}
\theta_2\!\!\!&=&\!\!\!\frac{3(u\vee v)+(u\wedge v)-2+\sqrt{(u+v)^2+4(1-u)(1-v)+2t}}{4},\\\label{theta3}
\theta_3\!\!\!&=&\!\!\!\frac{4(u\vee v)+2(u\wedge v)-3+\sqrt{16(u\vee v)^2+4(u\wedge v)^2-24(u\vee v)-12(u\wedge v)+16uv+7t+16}}{7},\\\label{theta4}
\theta_4\!\!\!&=&\!\!\!\frac{5(u\vee v)+3(u\wedge v)-4+\sqrt{4(u\vee v)^2+16(u\wedge v)^2-12(u\vee v)-24(u\wedge v)+16uv+7t+16}}{7},\\\label{theta5}
\theta_5\!\!\!&=&\!\!\!\frac{3(u+v-1)+\sqrt{3(5u^2+5v^2-6u-6v+2uv+2t+5)}}{6},
\end{eqnarray}
and $\chi_{\theta_i}:=\chi_{\theta_i}(u,v)$, $i=1,\ldots,5$, are the characteristic functions of the five respective subsets identified by the inequalities in \eqref{gammainf}.
\end{theorem}

\begin{proof}We determine the pointwise supremum. First, we compute the integral
\begin{equation}\label{Intunder}
I_1:=\int_0^1\underline{C}_{(a,b),\theta}(u,1-u)\,du.
\end{equation}
Note that we have to study three cases---numbered I,\, II and III---, depending on the location of the point $(a,b)$ in $[0,1]^2$. Figure \ref{fig:3under} shows these three cases, where the mass distribution of the copula $\underline{C}_{(a,b),\theta}$ is represented with continuous lines and the point $(a,b)$ with a dot, and Figure \ref{fig:caseI} shows the values of the copula $\underline{C}_{(a,b),\theta}$ in case I (the cases II and III can be described in a similar way and are omitted).
\begin{figure}[htb]
\begin{center}\epsfig{file=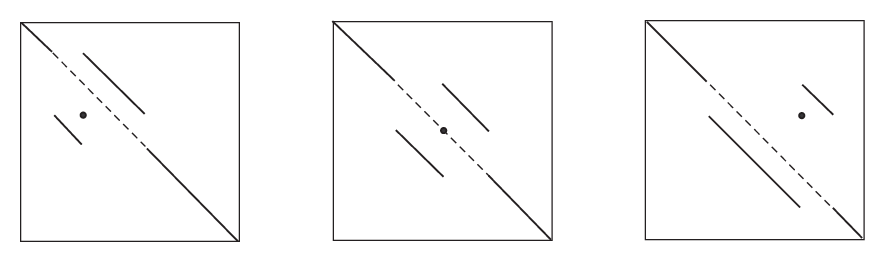,width=13cm}
\caption{The respective cases I, II and III (from left to right) for computing the integral \eqref{Intunder}.}\label{fig:3under}
\end{center}
\end{figure}
\begin{figure}[htb]
\begin{center}\epsfig{file=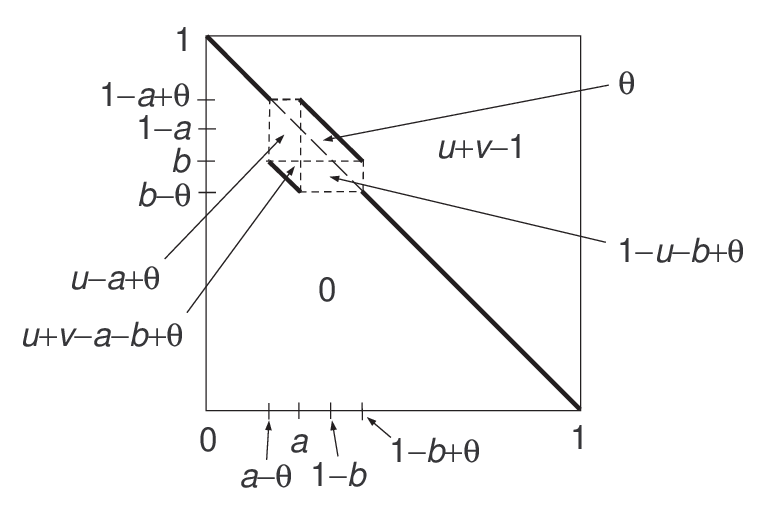,width=7cm}
\caption{The values of the copula $\underline{C}_{(a,b),\theta}$ in case I for computing \eqref{Intunder}.}\label{fig:caseI}
\end{center}
\end{figure}
\begin{list}{}{}
\item[Case I.] $I_1=\displaystyle\int_{a-\theta}^{a}(u-a+\theta)\,du+\int_{a}^{1-b}\theta\,du+\int_{1-b}^{1-b+\theta}(1-u-b+\theta)\,du$.
\item[Case II.] $I_1=\displaystyle\int_{a-\theta}^{a}(u-a+\theta)\,du+\int_{a}^{1-b+\theta}(1-u-b+\theta)\,du$.
\item[Case III.] $I_1=\displaystyle\int_{a-\theta}^{1-b}(u-a+\theta)\,du+\int_{1-b}^{a}(1-a-b+\theta)\,du+\int_{a}^{1-b+\theta}(1-u-b+\theta)\,du$.
\end{list}
In the three cases we obtain
$$I_1=\theta(1-a-b+\theta).$$

Now we compute the integral
\begin{equation}\label{Intunder2}
I_2:=\displaystyle\int_0^1\underline{C}_{(a,b),\theta}(u,u)\,du.
\end{equation}
We have to consider ten cases depending on the location of the point $(a,b)$, but, if we suppose $b\le a$, by symmetry, we only need to consider five. Figure \ref{fig:5casesinf} shows those five cases---numbered from 1 to 5---, and Figure \ref{fig:case2} shows the values of the copula $\underline{C}_{(a,b),\theta}$ in case 2 (the rest of the cases can be described in a similar way and are omitted).
\begin{figure}[htb]
\begin{center}\epsfig{file=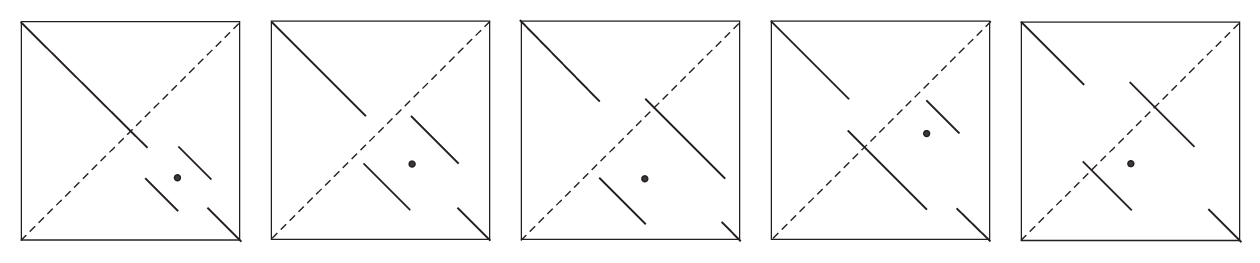,width=17.5cm}
\caption{The respective cases $1,2,3,4$ and $5$ (from left to right) for computing the integral in \eqref{Intunder2}.}\label{fig:5casesinf}
\end{center}
\end{figure}
\begin{figure}[htb]
\begin{center}\epsfig{file=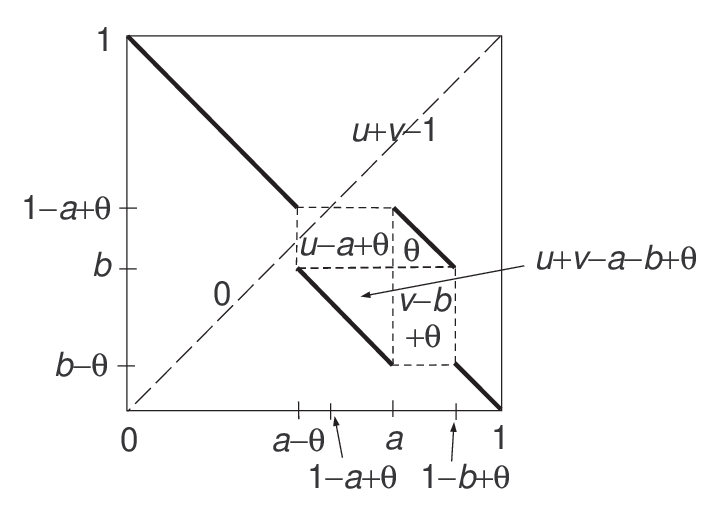,width=6.5cm}
\caption{The values of the copula $\underline{C}_{(a,b),\theta}$ in case 2 for computing \eqref{Intunder2}. Note that the inequalities $a-\theta\le 1-a+\theta$, $a-\theta\ge b$ and $1-a+\theta\le a$ must be satisfied in this case.}\label{fig:case2}
\end{center}
\end{figure}
Thus, computing the integral $I_2$ for these five cases we obtain:
\begin{list}{}{}
\item[Case 1.]
\begin{eqnarray*}
I_2&=&\displaystyle\int_{1/2}^1(2u-1)\,du\\
&=&\frac{1}{4}.
\end{eqnarray*}

\item[Case 2.]
\begin{eqnarray*}
I_2&=&\displaystyle\int_{a-\theta}^{1-a+\theta}(u-a+\theta)\,du+\int_{1-a+\theta}^1(2u-1)\,du\\
&=&(a-\theta)^2-(a-\theta)+\frac{1}{2}.
\end{eqnarray*}

\item[Case 3.]
\begin{eqnarray*}
I_2&=&\displaystyle\int_{a-\theta}^{a}(u-a+\theta)\,du+\int_a^{(1+\theta)/2}\theta\,du+\int_{(1+\theta)/2}^1(2u-1)\,du\\
&=&\frac{1}{4}+\frac{3}{4}\theta^2+\theta\left(\frac{1}{2}-a\right).
\end{eqnarray*}

\item[Case 4.]
\begin{eqnarray*}
I_2&=&\displaystyle\int_{(a+b-\theta)/2}^{b}(2u-a-b+\theta)\,du+\int_b^{1-a+\theta}(u-a+\theta)\,du+\int_{1-a+\theta}^1(2u-1)\,du\\
&=&\frac{3a^2+2a(b-3\theta-2)-b^2-2b\theta+3\theta^2+4\theta+2}{4}.
\end{eqnarray*}

\item[Case 5.]
\begin{eqnarray*}
I_2&=&\displaystyle\int_{(a+b-\theta)/2}^b(2u-a-b+\theta)\,du+\int_b^a(u-a+\theta)\,du
+\int_a^{(1+\theta)/2}\theta\,du+\int_{(1+\theta)/2}^{1}(2u-1)\,du\\
&=&\frac{1}{4}\left(1-(a-b)^2\right)+\frac{1}{2}(1-a-b)\theta+\frac{\theta^2}{2}.
\end{eqnarray*}
\end{list}

Thus we have
\begin{equation*}\label{Intinf}
I_2=\left\{\begin{array}{lllll}
\frac{1}{4}, & a\ge \frac{1}{2}+\theta, &\\ \noalign{\medskip}
\frac{1}{4}+\left(a-\theta-\frac{1}{2}\right)^2,& \left( b+\theta\right)\vee \frac{1+\theta}{2}\leq a\leq \frac{1}{2}+\theta, &\\ \noalign{\medskip}
\frac{1}{4}(1+2\theta-4a\theta+3\theta^2),& b+\theta \leq a \leq \frac{1+\theta}{2},&\\ \noalign{\medskip}
\frac{1}{4}((\theta+1-a-b)(3\theta-3a+b+1)+1),& \frac{1+\theta}{2} \leq a\leq  b+\theta,&\\ \noalign{\medskip}
\frac{1}{4}\left(1-(a-b)^2\right)+\frac{1}{2}(1-a-b)\theta+\frac{\theta^2}{2},& a \leq \left(b+\theta\right)\wedge \frac{1+\theta}{2}.&
\end{array} \right.
\end{equation*}

Hence, the values that $\gamma\left(\underline{C}_{(a,b),\theta}\right)$ can take are given by
\begin{equation}\label{gammainf}
\gamma\left(\underline{C}_{(a,b),\theta}\right)\!=\!\left\{\begin{array}{llll}
 4\theta^2+4\theta(1-a-b)-1, &  \frac{1}{2}+\theta \leq a\vee b,\\ \noalign{\medskip}
 \left(2(a\vee b)-2\theta-1\right)^2+4\theta^2+4\theta(1-a-b)-1, & \left((a\wedge b)+\theta\right)\vee \frac{1+\theta}{2} \leq a\vee b\leq \frac{1}{2}+\theta,\\ \noalign{\medskip}
 2\theta-4(a\vee b)\theta+7\theta^2+4\theta(1-a-b)-1, &  (a\wedge b)+\theta \leq a\vee b\leq \frac{1+\theta}{2},\\ \noalign{\medskip}
 (a+b-1-4\theta)^2+2(a+b-1-\theta)|a-b|-9\theta^2-1, &  \frac{1+\theta}{2}\leq a\vee b\le (a\wedge b)+\theta,\\ \noalign{\medskip}
 6\theta^2+6\theta(1-a-b)-(a-b)^2-1, &  a\vee b\leq \left((a\wedge b)+\theta\right)\wedge \frac{1+\theta}{2}.
\end{array} \right.
\end{equation}

Now, let $C$ be a copula such that $\gamma(C)=t$, with $t\in[-1,1]$. Let $\theta = C(a,b)$, where $(a,b)$ is fixed. Then ${\underline{C}}_{(a,b),\theta}\le C$. Since $\gamma\left(\underline{C}_{(a,b),\theta}\right)\le \gamma(C)$, then it holds $\gamma\left(\underline{C}_{(a,b),\theta}\right)\le t$. By considering the first value of $\gamma\left(\underline{C}_{(a,b),\theta}\right)=4\theta^2+4\theta(1-a-b)-1$ in \eqref{gammainf}, we have
$$4\theta^2+4\theta(1-a-b)-1\le t.$$
Let $g(\theta)=4\theta^2+4\theta(1-a-b)-1-t$. Then $\theta$ must be less or equal than the greatest of the roots of $g$, which is given by
$$\theta_1=\frac{a+b-1+\sqrt{(a+b-1)^2+t+1}}{2}$$
(note that by changing $(a,b)$ by $(u,v)$ we obtain \eqref{theta1}). A similar reasoning with the other values in \eqref{gammainf} leads us to $\theta_2, \theta_3,\theta_4$ and $\theta_5$ given by \eqref{theta2}, \eqref{theta3}, \eqref{theta4} and \eqref{theta5}, respectively. It follows that $\theta\le \max\left(\chi_{\theta_1}\theta_1,\chi_{\theta_2}\theta_2,\chi_{\theta_3}\theta_3,\chi_{\theta_4}\theta_4,\chi_{\theta_5}\theta_5\right)$, so that $C(a,b)\le\min\left(a,b,\max\left(\chi_{\theta_1}\theta_1,\chi_{\theta_2}\theta_2,\chi_{\theta_3}\theta_3,\chi_{\theta_4}\theta_4,\chi_{\theta_5}\theta_5\right)\right)$. \par

To establish \eqref{Gsup}, it is sufficient to show that for each pair $(a,b)$ there exists a copula $C$ in ${\bf G}_t$ such that
$$C(a,b)=\min\left(a,b,\max\left(\chi_{\theta_1}\theta_1,\chi_{\theta_2}\theta_2,\chi_{\theta_3}\theta_3,\chi_{\theta_4}\theta_4,\chi_{\theta_5}\theta_5\right)\right).$$
Assume $a\le b$ (the case $a>b$ is similar and we omit it), and choose $i\in\{1,2,3,4,5\}$ such that $\theta_i= \max\left(\theta_1,\theta_2,\theta_3,\theta_4,\theta_5\right)$. If $\theta_i\le a$, then $\underline{C}_{(a,b),\theta_i}\in {\bf G}_t$ and $\underline{C}_{(a,b),\theta_i}(a,b)=\theta_i$. If $\theta_i > a$, then $C_\alpha=\alpha M+(1-\alpha)\underline{C}_{(a,b),a}$, for $\alpha\in[0,1]$, is a family of copulas for which $C_\alpha (a,b) = a$ and $\gamma\left(C_\alpha\right)$ is a continuous function of $\alpha$ satisfying
$\gamma\left(C_0\right) = \gamma\left(\underline{C}_{(a,b),a}\right)\le t \le 1 = \gamma(M) = \gamma\left(C_1\right)$. It now follows from the intermediate
value theorem that there is $\alpha\in [0,1]$ such that $\gamma\left(C_\alpha\right)=t$.

For the pointwise infimum, note that the expression \eqref{Ginf} follows from Proposition \ref{relation}, which completes the proof.
\end{proof}

\begin{remark}Due to the complexity of its calculation and expression, the explicit formula for \eqref{Gsup} is provided in Appendix. Observe also that, as a consequence of Proposition \ref{relation}, we have that Equation \eqref{Ginf} can also be expressed as $\underline{G}_t(u,v)=u-\overline{G}_{-t}(u,1-v)$ for all $(u,v)\in[0,1]^2$.
\end{remark}

We know that the functions $\underline{G}_t$ and $\overline{G}_t$ given by \eqref{Ginf} and \eqref{Gsup}, respectively, are quasi-copulas, and unlike the pointwise infimum and supremum for the measures Kendall's tau, Spearman's rho and Blomqvist's beta, and the best-possible lower bound for the case of the Spearman's footrule, they can be proper quasi-copulas. We start studying $\overline{G}_t$.

\begin{proposition}\label{casesup}Let $\overline{G}_t$ be the quasi-copula given by \eqref{Gsup}. Then it holds that:
\begin{enumerate}
\item[{\rm (a)}] $\overline{G}_{-1}=W$.
\item[{\rm (b)}] $\overline{G}_t$ is a proper quasi-copula if, and only if, $-1<t<0$.
\item[{\rm (c)}] $\overline{G}_t$ is a copula, different from $M$ and $W$, if, and only if, $0\le t<1/2$.
\item[{\rm (d)}] $\overline{G}_t=M$ if, and only if, $1/2\le t\le 1$.
\end{enumerate}
\end{proposition}

\begin{proof}First, we prove part (a). When $t=-1$, we have that the regions ${\cal{R}}_2$ and ${\cal{R}}_5$ are reduced to segments (see Appendix). Assume $u\le v$, then we have regions ${\cal{R}}_1=[0,\frac{1}{2}]\times[\frac{1}{2},1]$, ${\cal{R}}_3\{(u,v)\in[0,1]^2: 0\le u\le 1/2,\, u\le v\le 1/2\}$ and ${\cal{R}}_4=\{(u,v)\in[0,1]^2: 1/2\le u\le 1,\, u\le v\le 1\}$. Therefore,
\begin{eqnarray*}
\min\left(u,\max\left(\frac{u+v-1+\sqrt{(u+v-1)^2}}{2}\right)\right)\!\!\!&=&\!\!\!\min\left(u,\max\left(\frac{u+v-1+|u+v-1|}{2}\right)\right)\\
&=&\!\!\!\max(0,u+v-1)\,\,\,{\rm in}\,\,\,{\cal{R}}_1,\\
\min\left(u,\max\left(\frac{2u+4v-3+\sqrt{(2u+4v-3)^2}}{7}\right)\right)\!\!\!&=&\!\!\!\min\left(u,\max\left(\frac{2u+4v-3+|2u+4v-3|}{7}\right)\right)\\
&=&\!\!\!0\,\,\,{\rm in}\,\,\,{\cal{R}}_3,\\
\min\left(u,\max\left(\frac{5v+3u-4+\sqrt{(4u+2v-3)^2}}{7}\right)\right)\!\!\!&=&\!\!\!\min\left(u,\max\left(\frac{5v+3u-4+|4u+2v-3|}{7}\right)\right)\\
&=&\!\!\!u+v-1\,\,\,{\rm in}\,\,\,{\cal{R}}_4,
\end{eqnarray*}
whence we easily obtain $\overline{G}_{-1}=W$.

Before proving part (b), we prove (c). For $0\le t<1/2$ we have
$$\overline{G}_t(u,v)=\min\left(u,v,\frac{3(u+v-1)+\sqrt{3(5u^2+5v^2-6u-6v+2uv+2t+5)}}{6}\right)$$
for all $(u,v)\in[0,1]^2$, since only region ${\cal{R}}_5$ exists (see Appendix). Thus,
$$\overline{G}_t(u,v)=\left\{\begin{array}{ll}
\displaystyle\frac{3(u+v-1)+\sqrt{3(5u^2+5v^2-6u-6v+2uv+2t+5)}}{6}, & (u+v)^2+2uv -6(u\wedge v)\le -1-t,\\ \noalign{\medskip}
M(u,v),& {\rm otherwise}.
\end{array} \right.$$
The set $S=\left\{(u,v)\in[0,1]^2 : (u+v)^2+2uv -6(u\wedge v)\le -1-t \right\}$ is the region between two symmetric hyperbolic arcs with respect to the line $v=u$, and cut at the points
$$p_1=\left(\frac{3+\sqrt{3-6t}}{6},\frac{3+\sqrt{3-6t}}{6}\right)\quad{\rm and}\quad p_2=\left(\frac{3-\sqrt{3-6t}}{6},\frac{3-\sqrt{3-6t}}{6}\right).$$
Therefore, the only $\overline{G}_t-$volumes to be studied are those of the rectangles that have some vertex in the interior of $S$---denoted by ${\rm int}(S)$---, and in turn, the study of these volumes can be reduced to the case in which the four vertices of the rectangle are in ${\rm int}(S)$. But these volumes are non-negative since, if $(u,v)\in {\rm int}(S)$, then we have
$$\frac{\partial^2 \overline{G}_t}{\partial u\partial v}(u,v)=\frac{\sqrt{3(t-12uv+6u+6v-2)}}{3\left(5u^2+5v^2-6u-6v+2uv+2t+5\right)^{3/2}}\ge 0,$$
i.e., $\underline{G}_t$ is 2-increasing at such points, and hence $\overline{G}_t$ is a copula, which proves part (c).

To prove that $\overline{G}_t$ is a proper quasi-copula for $-1<t<0$---region ${\cal{R}}_1$ does not exist for $t>-\frac{3}{4}$, region ${\cal{R}}_2$ does not exist for $t>-\frac{4}{9}$, and regions ${\cal{R}}_3$ and ${\cal{R}}_4$ do not exist for $t>-\frac{4}{13}$ (see Appendix)---, note that the point $p_1$ (similarly, the point $p_2$) belongs to the boundary of the set $S$, and we have
$$\frac{\partial^2 \overline{G}_t}{\partial u\partial v}\left(\frac{3+\sqrt{3-6t}}{6},\frac{3+\sqrt{3-6t}}{6}\right)=\frac{t}{3}<0.$$

Finally, since all the regions ${\cal{R}}_i$, $i=1,\ldots,5$, do not exist for $t\in[\frac{1}{2},1]$ (see Appendix), we have $\overline{G}_t=M$, which completes the proof.
\end{proof}

\begin{remark}
We note that the copulas in Proposition \ref{casesup}(c) has a singular component (\cite{Durante2016book,Nelsen2006}) on two segments on the main diagonal and on the two hyperbolic arcs of the set $S$ between the points $p_1$ and $p_2$. Figure \ref{fig:G} shows the graph and the level curves of the copula $\overline{G}_0$.
\begin{figure}[htb]
\begin{center}\epsfig{file=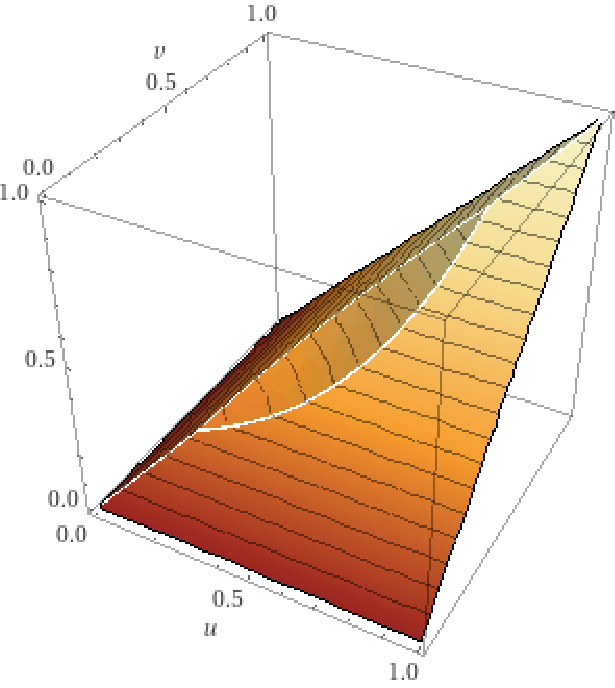,width=5cm}\hspace{2cm}\epsfig{file=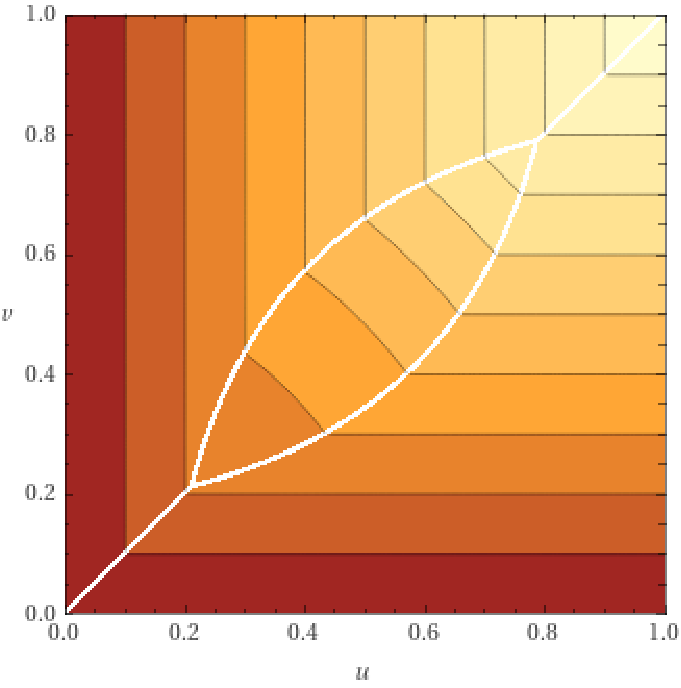,width=5cm}
\caption{The graph (left) and the level curves (right) of the copula $\overline{G}_0$.}\label{fig:G}
\end{center}
\end{figure}
\end{remark}

The next result show the study of the lower bound, whose proof is immediate from Equation \eqref{Ginf} and Proposition \ref{casesup}.

\begin{proposition}Let $\underline{G}_t$ be the quasi-copula given by \eqref{Ginf}. Then it holds that:
\begin{enumerate}
\item[{\rm (a)}] $\underline{G}_{1}=M$.
\item[{\rm (b)}] $\underline{G}_t$ is a proper quasi-copula if, and only if, $0<t<1$.
\item[{\rm (c)}] $\underline{G}_t$ is a copula, different from $M$ and $W$, if, and only if, $-1/2< t\le 0$.
\item[{\rm (d)}] $\underline{G}_t=W$ if, and only if, $-1\le t\le -1/2$.
\end{enumerate}
\end{proposition}

\section{Conclusions}\label{conclus}

In this note, we have found best-possible bounds ---supremum ($\overline{G}_{t}$) and infimum ($\underline{G}_{t}$)--- on the set of copulas with a given value of the well-known measure of association Gini's gamma. We showed that, unlike the cases for the measures Kendall's tau, Spearman's rho and Blomqvist's beta, and the lower bound for the case of the Spearman's footrule, the pointwise supremum (respectively, infimum) $\overline{G}_t$ (respectively, $\underline{G}_t$) can be a proper quasi-copula for $-1 < t< 0$ (respectively, $0<t<1$), being copulas for the rest of the cases.
\bigskip\medskip

\noindent{\large\bf Acknowledgement}
\medskip

The author acknowledges the support of the program FEDER-Andaluc\'ia 2014-2020 under research project UAL2020-AGR-B1783 and project PID2021-122657OB-I00 by the Ministerio de Ciencia e Innovaci\'on (Spain).

\section*{Appendix}

In order to provide an explicit formula of expression \eqref{Gsup}, we expand $\chi_{\theta_i}\theta_i$, for $i=1,\ldots,5$. We show the case $\chi_{\theta_1}\theta_1$ in detail, following the rest of the cases a similar reasoning. From \eqref{gammainf} we have
\begin{equation}\label{eq:case1}
\gamma\left(\underline{C}_{(a,b),\theta}\right)= 4\theta^2+4\theta(1-a-b)-1\quad{\rm if}\quad \frac{1}{2}+\theta \leq a\vee b.
\end{equation}
Assume $a\le b$ (the case $a>b$ is similar and we omit it). Since $\underline{C}_{(a,b),\theta}(a,b)=\theta$ and $4\theta^2+4\theta(1-a-b)-1$ in \eqref{eq:case1} is increasing in $\theta$, the maximum of $\gamma\left(\underline{C}_{(a,b),\theta}\right)$ is reached at $\theta=a$, whence $4\theta^2+4\theta(1-a-b)-1=4a(1-b)-1\ge t$ implies
\begin{equation}\label{eq:b1}
b\le 1-\frac{1+t}{4a}.
\end{equation}
Now, solving the equation $4\theta^2+4\theta(1-a-b)-1=t$ in $\theta$ we obtain the greatest root
$$\theta=\frac{a+b-1+\sqrt{(a+b-1)^2+t+1}}{2}.$$
Since
$\frac{1}{2}+\theta \leq b$---recall \eqref{eq:case1}---then
$$\frac{1}{2}+\frac{a+b-1+\sqrt{(a+b-1)^2+t+1}}{2}\le b$$
implies
\begin{equation}\label{eq:b2}
b\ge \frac{2a-t-2}{4a-2}
\end{equation}
as long as $a\le \frac{1}{2}$. Following \eqref{eq:b1} and \eqref{eq:b2}, taking into account the case $a>b$ as well, and renaming $u:=a$ and $v:=b$, we can conclude that
$$\chi_{\theta_1}\theta_1=\frac{u+v-1+\sqrt{(u+v-1)^2+t+1}}{2}\quad{\rm if}\quad (u,v)\in{\mathcal{R}}_1,$$
where
$${\mathcal{R}}_1=\left\{(u,v)\in[0,1]^2:\frac{2(u\wedge v)-t-2}{4(u\wedge v)-2}\le (u\vee v)\le 1-\frac{1+t}{4(u\wedge v)},\,(u\wedge v)\le \frac{1}{2}\right\}.$$
Finally, we note that the set ${\mathcal{R}}_1$ exists---i.e., it is non-empty---for all $t\in\left[-1,-\frac{3}{4}\right]$ since, e.g., solving the inequality
$$\frac{2u-t-2}{4u-2}\le 1-\frac{1+t}{4u}$$
in $u$ we obtain
$$\left(u\le \frac{1-\sqrt{-4t-3}}{4}\,\,{\rm and}\,\, u\le\frac{1+\sqrt{-4t-3}}{4}\right)\,\,{\rm or}\,\,\left(u\le \frac{1-\sqrt{-4t-3}}{4}\,\,{\rm and}\,\, u\le\frac{1+\sqrt{-4t-3}}{4}\right)$$
(recall $u\le \frac{1}{2}$).

Similarly, we expand $\chi_{\theta_i}\theta_i$, for $i=2,\ldots,5$:
$$\chi_{\theta_2}\theta_2=\frac{3(u\vee v)+(u\wedge v)-2+\sqrt{(u+v)^2+4(1-u)(1-v)+2t}}{4}\quad{\rm if}\quad (u,v)\in{\mathcal{R}}_2,$$
where
\begin{eqnarray*}
{\mathcal{R}}_2\!\!\!&=&\!\!\!\left\{(u,v)\in[0,1]^2:\left(\frac{12(u\wedge v)^2-8(u\wedge v)-t}{8(u\wedge v)-4}\vee\frac{2(u\wedge v)+2+\sqrt{4(u\wedge v)^2-4(u\wedge v)+3t+4}}{6}\right)\right.\\
&&\left.\le (u\vee v)\le\left(\frac{2(u\wedge v)-t-2}{4(u\wedge v)-2}\wedge \frac{3(u\wedge v)+1-\sqrt{5(u\wedge v)^2-2(u\wedge v)+t+1}}{2}\right)\right\}.
\end{eqnarray*}
The set ${\mathcal{R}}_2$ exists for all $t\in[-1,-\frac{4}{9}]$.
\begin{eqnarray*}
\chi_{\theta_3}\theta_3\!\!\!&=&\!\!\!\frac{4(u\vee v)+2(u\wedge v)-3+\sqrt{16(u\vee v)^2+4(u\wedge v)^2-24(u\vee v)-12(u\wedge v)+16uv+7t+16}}{7}\\
&&{\rm if}\quad (u,v)\in{\mathcal{R}}_3,
\end{eqnarray*}
where
\begin{eqnarray*}
{\mathcal{R}}_3\!\!\!&=&\!\!\!\left\{(u,v)\in[0,1]^2:3-5(u\wedge v)-\sqrt{36(u\wedge v)^2-36(u\wedge v)-t+8}\le (u\vee v)\right.\\
&&\left.\le \left(\frac{3(u\wedge v)^2+6(u\wedge v)-t-1}{8(u\wedge v)}\wedge \frac{2(u\wedge v)+2+\sqrt{4(u\wedge v)^2-4(u\wedge v)+3t+4}}{6}
\right)\right\}.
\end{eqnarray*}
The set ${\mathcal{R}}_3$ exists for all $t\in[-1,-\frac{4}{13}]$.
\begin{eqnarray*}
\chi_{\theta_4}\theta_4\!\!\!&=&\!\!\!\frac{5(u\vee v)+3(u\wedge v)-4+\sqrt{4(u\vee v)^2+16(u\wedge v)^2-12(u\vee v)-24(u\wedge v)+16uv+7t+16}}{7}\\
&&{\rm if}\quad (u,v)\in{\mathcal{R}}_4,
\end{eqnarray*}
where
\begin{eqnarray*}
{\mathcal{R}}_4\!\!\!&=&\!\!\!\left\{(u,v)\in[0,1]^2:\frac{5(u\wedge v)+3+\sqrt{9(2(u\wedge v)-1)^2+11(t+1)}}{11}\le (u\vee v)\right.\\
&&\left.\le \left(\frac{12(u\wedge v)^2-8(u\wedge v)-t}{8(u\wedge v)-4}\wedge \frac{4(u\wedge v)+2-\sqrt{16(u\wedge v)^2-8(u\wedge v)+3t+4}}{3}\right)\right\}.
\end{eqnarray*}
The set ${\mathcal{R}}_4$ exists for all $t\in[-1,-\frac{4}{13}]$.
\begin{eqnarray*}
\chi_{\theta_5}\theta_5=\frac{3(u+v-1)+\sqrt{3(5u^2+5v^2-6u-6v+2uv+2t+5)}}{6}\quad{\rm if}\quad (u,v)\in{\mathcal{R}}_5,
\end{eqnarray*}
where
\begin{eqnarray*}
{\mathcal{R}}_5\!\!\!&=&\!\!\!\left\{(u,v)\in[0,1]^2:(u\vee v)\le \left(3-5(u\wedge v)-\sqrt{36(u\wedge v)^2-36(u\wedge v)-t+8}\right.\right.\\
&&\left.\left.\wedge \frac{5(u\wedge v)+3+\sqrt{9(2(u\wedge v)-1)^2+11(t+1)}}{11}\wedge -2(u\wedge v)+\sqrt{3(u\wedge v)((u\wedge v)+2)-(t+1)}\right)\right\}
\end{eqnarray*}
The set ${\mathcal{R}}_5$ exists for all $t\in[-1,\frac{1}{2}]$.

\end{document}